\documentclass{aptpub}

\authornames{S.B. CONNOR} 
\shorttitle{Optimal co-adapted coupling on $K_{\lowercase{n}}^{\lowercase{d}}$} 

\usepackage{tikz}
\usepackage{xspace}
\usepackage{fancybox}

\newcommand{\abs}[1]{\left\vert#1\right\vert}
\newcommand{\set}[1]{\left\{#1\right\}}

\newcommand{\squ}[1]{\left[#1\right]}
\newcommand{\bra}[1]{\left(#1\right)}
\newcommand{\norm}[1]{\left\Vert#1\right\Vert}
\renewcommand{\prob}[1]{\operatorname{\mathbb{P}}\left(#1\right)}

\newcommand{\Ex}[1]{\operatorname{\mathbb{E}}\left[#1\right]}

\renewcommand{\d}{\operatorname{d}\!}
\newcommand{\indev}[1]{\mathbf{1}_{\squ{#1}}}

\newcommand{\hyp}{\ensuremath{\mathbb{Z}_2^d}\xspace}
\newcommand{\law}{\ensuremath{\mathcal{L}}}
\newcommand{\RR}{\ensuremath{\mathbb R}}
\newcommand{\distname}[1]{\mbox{\textup{#1}}}
\newcommand{\expdist}{\distname{Exp}}
\newcommand{\NN}{\mathbb N}
\newcommand{\unifdist}{\distname{Uniform}}
\newcommand{\C}[1]{\mbox{[C#1]}}
\newcommand{\tmi}{{t-}}
\newcommand{\tm}{TM}
\newcommand{\hypK}{\ensuremath{K_n^d}\xspace}

\newcommand{\lc}[2][]{\ensuremath{\hat \lambda_t(#1 #2)}}
\newcommand{\rinf}{\rightarrow\infty}

\newcommand{\ie}{\emph{i.e.}\ }

\numberwithin{equation}{section}  

\begin{document}

\title{Optimal co-adapted coupling for a random walk on the hyper-complete-graph} 

\authorone[University of York]{Stephen Connor} 

\addressone{Department of Mathematics, University of York, York, YO10 5DD. UK} 

\begin{abstract}
The problem of constructing an optimal co-adapted coupling for a pair of symmetric random walks on \hyp was considered by Connor and Jacka \cite{Connor.Jacka-2008}, and the existence of a coupling which is stochastically fastest in the class of all such co-adapted couplings was demonstrated. In this paper we show how to generalise this construction to an optimal co-adapted coupling for the continuous-time symmetric random walk on \hypK, where $K_n$ is the complete graph with $n$ vertices. Moreover, we show that although this coupling is not maximal for any $n$ (\emph{i.e.} it does not achieve equality in the coupling inequality), it does tend to a maximal coupling as $n\rinf$. 
\end{abstract}




\section{Introduction}\label{sec:intro}

The concept of coupling two copies of a Markov chain in order to prove ergodicity statements dates back to \cite{Doeblin1938}, and is now a well-used and elegant technique. Suppose that we have a Markov process $X$ on some state space $\mathcal S$, and let $Y$ be a copy of $X$ started from a different initial distribution: a \emph{coupling} of these two processes is defined as follows.
\begin{defn}[Coupling]\label{defn:coupling}
        A \emph{coupling} of $X$ and $Y$ is a process $(X^c,Y^c)$ on $\mathcal S \times \mathcal S$ such that
            \[ X^c\stackrel{\mathcal{D}}{=} X \quad \text{and} \quad Y^c\stackrel{\mathcal{D}}{=} Y \,, \]
where $\stackrel{\mathcal{D}}{=}$ denotes equality in distribution.
      \end{defn}
That is, viewed marginally, $X^c$ behaves as a version of $X$, and $Y^c$ as a version of $Y$. The \emph{coupling time} $\tau^c$ is defined by
\[ \tau^c =  \inf\set{t\,:\, X^c_s = Y^c_s \text{ for all $s\geq t$}}\,, \]
and the coupling is called \emph{successful} if $\prob{\tau^c<\infty} = 1$.

Recall that the coupling inequality bounds the tail distribution of \emph{any} coupling of $X$ and $Y$ by the total variation distance between the two processes (see \emph{e.g.} \cite{Lindvall-2002}):
\begin{equation}\label{eqn:coupling-inequality}
\norm{\law(X_t)-\law(Y_t)}_{\textup{TV}} \leq \prob{\tau^c>t} \,, 
\end{equation}
where $\law(X_t)$ is the law of $X_t$. Moreover, there always exists a \emph{maximal} coupling of $X$ and $Y$; that is, one which achieves equality in \eqref{eqn:coupling-inequality} \cite{Griffeath-1975,Sverchkov.Smirnov-1990}. Thus, in order to obtain a good estimate of the rate at which the distributions of $X_t$ and $Y_t$ converge, it suffices to find a `good' coupling -- one which has a small coupling time $\tau^c$ -- and preferably one which is maximal. However, maximal couplings are, in most cases, unintuitive non-Markovian affairs, and extremely difficult to work with. It is therefore natural to consider the class of \emph{co-adapted} couplings. These may not be maximal, but are generally more intuitive; nearly all couplings used in practice fall into this class.





\begin{defn}[Co-adapted coupling]\label{defn:co-adapted}
  A coupling $(X^c,Y^c)$ is called \emph{co-adapted} if there exists a filtration $\bra{\mathcal{F}_t}_{t\geq 0}$ such that
        \begin{enumerate}
            \item $X^c$ and $Y^c$ are both adapted to $\bra{\mathcal{F}_t}_{t\geq 0}$\,;
            \item for any $0\leq s\leq t$,
                \[ \law\bra{X^c_t \,|\, \mathcal{F}_s} = \law\bra{X^c_t \,|\, X^c_s} \quad\text{and}\quad
                    \law\bra{Y^c_t \,|\, \mathcal{F}_s} = \law\bra{Y^c_t \,|\, Y^c_s}\,.
                \]
        \end{enumerate}
\end{defn}
In other words, $(X^c,Y^c)$ is co-adapted if $X^c$ and $Y^c$ are both Markov with respect to a common filtration: for co-adapted couplings it is always possible to make $X^c$ and $Y^c$ coalesce at the first collision time of the two chains. (Note that some authors (\emph{e.g.} \cite{Hsu.Sturm-,Kuwada-2007,Kuwada2009a}) use the term `Markovian' in place of `co-adapted', but this seems confusing when the \emph{joint} process $(X^c,Y^c)$ is not required to be Markov; we therefore prefer to reserve `Markovian' for couplings where this stronger condition is satisfied.)




With some families of chains it is possible to produce co-adapted couplings with coupling times of the same order of magnitude as the mixing time. For example, for simple random walk on $\hyp$ the mixing time is $\log d/4$ and the coupling of \cite{Connor.Jacka-2008} takes time $\log d/2$. However, as pointed out by \cite{Griffeath-1975}, in general a Markovian maximal coupling need not exist, and there are well-known chains for which co-adapted couplings perform significantly worse than maximal. For example, the transposition shuffle on the symmetric group $S_n$ and the Gibbs sampler on the $N$-simplex both have mixing times of order $O(n\log n)$, but it is impossible for co-adapted couplings to do better than $O(n^2)$ and $O(n^2\log n)$ respectively. It is therefore necessary in these instances to turn to non-co-adapted couplings in order to obtain a good bound on the mixing time -- see the preprints by \cite{Burton2011} and \cite{Smith2011} for details.

Due to the fact that maximal couplings are usually unintuitive and impractical to work with, whereas co-adapted couplings are used extensively both in theoretical problems  (\emph{e.g.} bounding mixing times) and in practical applications (\emph{e.g.} perfect simulation techniques), it is clearly of interest to obtain a better understanding of how good co-adapted couplings can be for various processes. This is a difficult question to answer in general, but explicit answers can be obtained in relatively simple cases: below we introduce one such process for which exact calculations are possible.


\subsection{Random walk on \hypK}\label{ssec:RW_on_complete}
Let $K_n$ be the complete graph with $n$ vertices (labelled $0,1,\dots,n-1$) and let \hypK be the set of $d$-tuples of the form $(x(1),\dots,x(d))$ with $x(i)\in K_n$, $1\leq i\leq d$. $\hypK$ forms a group under coordinate-wise addition modulo $n$, and $K_2^d \equiv \hyp$. A simple symmetric continuous-time random walk $X$ on \hypK may be defined by moving the $i^{th}$ coordinate of $X$ to a different, uniformly chosen value at incident times of a unit rate Poisson process, independently of all other coordinates. The unique equilibrium distribution of $X$ is the uniform distribution on \hypK. 

Suppose now that we have a coupled pair $(X^c,Y^c)$, and let
    \[ U^c_t = \set{1\leq i \leq d \,:\, X^c_t(i) \neq Y^c_t(i)}\quad\text{and}\quad  M^c_t = \set{1\leq i \leq d \,:\, X^c_t(i) = Y^c_t(i)} \]
respectively denote the sets of unmatched and matched coordinates at time $t\geq 0$. The coupling time $\tau^c$ clearly satisfies
    \[ \tau^c = \inf\set{t\geq 0 \,:\, X^c_s = Y^c_s \;\; \forall\, s\geq t} = \inf\set{t\geq 0\,:\, U_s^c = \emptyset \;\; \forall \, s\geq t} \,. \]



In \cite{Connor.Jacka-2008}, an explicit, intuitive coupling strategy is described when $n=2$, and is shown to yield the stochastically minimal coupling time of all co-adapted couplings. This coupling strategy at time $t$ depends only on the parity of $N_t=\abs{U_t}$, and may be summarised as follows: 
\begin{itemize}
\item matched coordinates are always made to move synchronously; \smallskip
\item if $N$ is odd, all unmatched coordinates of $X$ and $Y$ are made to evolve \emph{independently} until $N$ becomes even; \smallskip
\item if $N$ is even, unmatched coordinates are coupled in \emph{pairs} -- when an unmatched coordinate on $X$ flips (thereby making a new match), a different unmatched coordinate on $Y$ is flipped at the same instant (making a total of two new matches).
\end{itemize}

This motivates the following question: what is the optimal co-adapted coupling when $n>2$? Intuitively, we expect the optimal strategy of that paper to become inefficient as $n$ gets large, since the rate at which unmatched coordinates can be made to agree using either `independent' or  `pairwise' coupling (as described above) is proportional to $N/n$. In Section~\ref{sec:co-adapted} we show how to describe the problem of finding an optimal co-adapted coupling as an exercise in stochastic control, and solve this problem to once again obtain a stochastically minimal coupling time (in which new matches are made asymptotically at rate $N$ as $n\rinf$); the proof is deferred to the appendix.
In Section~\ref{sec:limits} we study the behaviour of this coupling as $d\to\infty$, for fixed $n$: we show that the time to stationarity exhibits a cutoff phenomenon which occurs strictly earlier than the coupling time, showing that our coupling is not maximal. As $n\to\infty$ however (with $d$ fixed), the optimal co-adapted coupling does tend to a maximal coupling. These are the results that we find most interesting since, to the best of the author's knowledge, this is the first time that such behaviour (a natural sequence of optimal co-adapted but sub-maximal couplings tending to a maximal coupling) has been observed. We conclude with some comments in Section~\ref{sec:comments} about a faster coupling for the random walks on \hypK, and about optimal couplings for random walks on $\mathbb{Z}_n^d$.


\section{Optimal co-adapted coupling on \hypK}\label{sec:co-adapted}

In order to find the optimal co-adapted coupling of $X$ and $Y$, it is first necessary to be able to describe a general coupling strategy $c\in\mathcal{C}$. To begin, let $\Lambda_{(i,k)(j,l)}$ ($1\leq i,j \leq d$ and $0\leq k,l \leq n-1$) be independent Poisson processes on $[0,\infty)$, each of rate $(n-1)^{-1}$. Now let $\set{W_{(i,k)(j,l)}}$ be a collection of piecewise constant processes on $[0,1]$, where $W_{(i,k)(j,l)}$ jumps at event times of $\Lambda_{(i,k)(j,l)}$ to values which are i.i.d. $\unifdist[0,1]$ (independently of all other $\Lambda$ and $W$ processes). We let $\bra{\mathcal{F}_t}_{t\geq 0}$ be any filtration to which all of the processes $\set{\Lambda_{(i,k)(j,l)},W_{(i,k)(j,l)}}$ are adapted.

The transitions of $X^c$ and $Y^c$ will be driven by the $W$ processes, and controlled by a process $\set{Q^c(t)}_{t\geq 0}$ which is adapted to $\bra{\mathcal{F}_t}_{t\geq 0}$, where
\[ Q^c(t)=\set{q^c_{(r,s)}(t)\,:\, 0\leq r,s\leq nd-1} \] is a $(nd)\times (nd)$ doubly-stochastic matrix.

A similar argument to that of \cite{Connor.Jacka-2008} shows that a general co-adapted coupling for $X$ and $Y$ may be defined as follows: if there is a jump in the process $W_{(i,k)(j,l)}$ at time $t\geq 0$, \emph{and} the value of $W_{(i,k)(j,l)}(t)$ satisfies $W_{(i,k)(j,l)}(t)\leq q^c_{([i-1]n+k,[j-1]n+l)}(t)$, then set $X^c_t(i)=k$ and $Y^c_t(j)=l$. To ease notation, in the sequel we shall write $q^c_{(i,k)(j,l)}(t)$ instead of $q^c_{([i-1]n+k,[j-1]n+l)}(t)$: thus $q^c_{(i,k)(j,l)}(t)$ is proportional to the instantaneous rate at which $(X^c_t(i),Y^c_t(j))$ jumps to $(k,l)$. (Note that this construction allows for the possibility of only one of $X^c$ and $Y^c$ actually changing its value at any given instant.)

Using the above construction, the rate at which $X^c(i)$ jumps from $r$ to $s\neq r$  is equal to
\[  \frac{1}{n-1} \,\sum_{j=1}^d \sum_{l=0}^{n-1} q^c_{(i,s)(j,l)}(t) = \frac{1}{n-1} \,, \]
since the double sum is simply the sum of the $([i-1]n+s)^{th}$ row of $Q^c(t)$, and hence equal to one. From this it follows directly that $X^c$ and $Y^c$ both
have the correct marginal transition rates to be continuous-time
simple random walks on \hypK as described in Section~\ref{ssec:RW_on_complete}, and are co-adapted.


\subsection{Stochastically optimal coupling}\label{ssec:optimal}

Our proposed optimal coupling $\hat c =\hat{c}_{n,d}$ once again depends upon the parity of $\hat N_t = |\hat U_t|$, the number of unmatched coordinates of $(\hat X, \hat Y) = (X^{\hat c},Y^{\hat c})$ at time $t$. It now also depends upon how this number relates to the parameter $n$.

\begin{defn}
\label{def:Q-hat}
The matrix process $\hat{Q}$ corresponding to the coupling $\hat{c}$ has non-zero entries given by the following rules.
\begin{itemize}
\item[\C1] $\hat{q}_{(i,k)(i,k)}(t) = 1$ for all $i\in \hat M_\tmi$ and all $k=0,\dots,n-1$; \smallskip
\item[\C2] if $\hat N_\tmi$ is even, \emph{or} $\hat N_\tmi\geq 2(n-1)/(n-2)$: for $i,j\in \hat U_\tmi$, with $i\neq j$, and for all $k,l\in\set{0,\dots,n-1}$,
  \begin{enumerate}
    \item[(i)] $\displaystyle{\hat{q}_{(i,k)(j,l)}(t) = \bra{\frac{1}{\hat N_\tmi-1}} \indev{k=Y_{t-}(i), \, l=X_{t-}(j)}}$;
    \item[(ii)] $\hat{q}_{(i,k)(i,k)}(t) = \indev{k\neq X_{t-}(i),\, k\neq Y_{t-}(i)}$;
    \end{enumerate} \smallskip
  \item[\C3] if $\hat N_\tmi$ is odd \emph{and} $\hat N_\tmi< 2(n-1)/(n-2)$, then $\hat{q}_{(i,k)(i,k)}(t) = 1$ for all $i\in \hat U_\tmi$ and all $k=0,\dots,n-1$.
   \end{itemize}

\end{defn}

Part \C1 of this definition ensures that no matches are ever broken under $\hat{c}$. The final two items define the strategy for making new matches. If $\hat N_\tmi$ is even, or else sufficiently large, we will see that \C2(i) implies that the rate at which two new matches are made is maximised; \C2(ii) then maximises the rate at which a single new match is made, subject to the constraint imposed by \C2(i). Finally, \C3 implies that if $\hat N_\tmi$ is odd, with $\hat N_\tmi<2(n-1)/(n-2)$, the coupling maximises the rate at which single matches are made. (Note that if $n=2$, \C3 applies whenever $\hat N_\tmi$ is odd; if $n=3$ then it applies when $\hat N_\tmi\in\set{1,3}$; while if $n\geq 4$, \C3 applies only when $\hat N_\tmi=1$.)

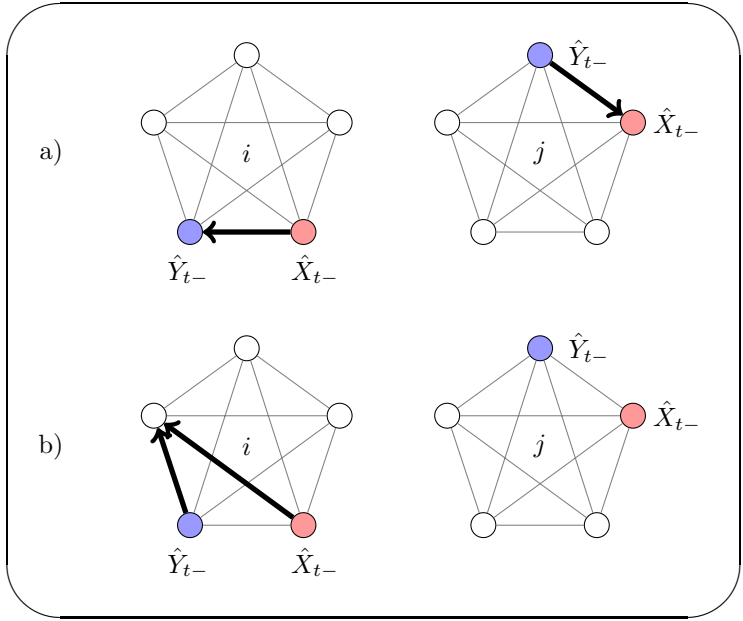
\begin{figure}[t]
\begin{center}
\ovalbox
{
\begin{tikzpicture}[scale=1.3]
\path (-2,0) node(a) {a)};
\path ( 0,1) node(1) [shape=circle,draw] {}
      ( -0.95,0.31) node(5) [shape=circle,draw] {}
      ( -0.58,-0.81) node(4) [shape=circle,draw,fill=blue!40] {}
      ( 0.58,-0.81) node(3) [shape=circle,draw,fill=red!40] {}
      (0.95,0.31) node(2) [shape=circle,draw] {};
\draw[color=gray] (1) -- (2) -- (3) -- (4) -- (5) -- (1) -- (3) -- (5) -- (2) -- (4) -- (1);

\path (0.7,-1.2) node(x1) [shape=circle] {$\hat X_{t-}$}
      (-0.6,-1.2) node(y1) [shape=circle] {$\hat Y_{t-}$}
      (0,0) node(i) [shape=circle] {$i$};
\draw[->,line width=2pt] (3) -- (4);

\path ( 3,1) node(1a) [shape=circle,draw,fill=blue!40] {}
      ( 2.05,0.31) node(5a) [shape=circle,draw] {}
      ( 2.42,-0.81) node(4a) [shape=circle,draw] {}
      ( 3.58,-0.81) node(3a) [shape=circle,draw] {}
      (3.95,0.31) node(2a) [shape=circle,draw,fill=red!40] {};
\draw[color=gray] (1a) -- (2a) -- (3a) -- (4a) -- (5a) -- (1a) -- (3a) -- (5a) -- (2a) -- (4a) -- (1a);

\path (4.4,0.3) node(x2) [shape=circle] {$\hat X_{t-}$}
      (3.5,1) node(y2) [shape=circle] {$\hat Y_{t-}$}
      (3,0) node(j) [shape=circle] {$j$};
\draw[->,line width=2pt] (1a) -- (2a);

\path (-2,-3) node(b) {b)};
\path ( 0,-2) node(1b) [shape=circle,draw] {}
      ( -0.95,-2.69) node(5b) [shape=circle,draw] {}
      ( -0.58,-3.81) node(4b) [shape=circle,draw,fill=blue!40] {}
      ( 0.58,-3.81) node(3b) [shape=circle,draw,fill=red!40] {}
      (0.95,-2.69) node(2b) [shape=circle,draw] {};
\draw[color=gray] (1b) -- (2b) -- (3b) -- (4b) -- (5b) -- (1b) -- (3b) -- (5b) -- (2b) -- (4b) -- (1b);

\path (0.7,-4.2) node(x1b) [shape=circle] {$\hat X_{t-}$}
      (-0.6,-4.2) node(y1b) [shape=circle] {$\hat Y_{t-}$}
      (0,-3) node(ib) [shape=circle] {$i$};
\draw[->,line width=2pt] (3b) -- (5b);
\draw[->,line width=2pt] (4b) -- (5b);

\path ( 3,-2) node(1ab) [shape=circle,draw,fill=blue!40] {}
      ( 2.05,-2.69) node(5ab) [shape=circle,draw] {}
      ( 2.42,-3.81) node(4ab) [shape=circle,draw] {}
      ( 3.58,-3.81) node(3ab) [shape=circle,draw] {}
      (3.95,-2.69) node(2ab) [shape=circle,draw,fill=red!40] {};
\draw[color=gray] (1ab) -- (2ab) -- (3ab) -- (4ab) -- (5ab) -- (1ab) -- (3ab) -- (5ab) -- (2ab) -- (4ab) -- (1ab);

\path (4.4,-2.7) node(x2b) [shape=circle] {$\hat X_{t-}$}
      (3.5,-2) node(y2b) [shape=circle] {$\hat Y_{t-}$}
      (3,-3) node(jb) [shape=circle] {$j$};
\end{tikzpicture}
}
\end{center}
\caption{Optimal co-adapted coupling strategy $\hat{c}_n$, when $n=5$: in both diagrams, the values of $\hat X_{t-}$ and $\hat Y_{t-}$ on coordinates $i,j\in \hat U_{t-}$ are highlighted. Suppose that $\hat X(i)$ jumps to a new value at time $t$. In a), this new value happens to equal $\hat Y_{t-}(i)$ -- in this case we select another unmatched coordinate ($j$) and move $\hat Y(j)$ to the value of $\hat X_{t-}(j)$, making a total of {\bf two} new matches. In b), the value of $\hat X_t(i)$ does not agree with that of $\hat Y_{t-}(i)$ -- here $\hat Y(i)$ is also made to jump so that $\hat X_t(i) = \hat Y_t(i)$, making {\bf one} new match.}\label{fig:intuition}
\end{figure}

Informally, when $n\geq 4$ and $\hat N_\tmi\geq 2$, $\hat{c}$ couples $\hat X$ and $\hat Y$ as follows (see Figure~\ref{fig:intuition}). If an unmatched coordinate $\hat X(i)$ jumps to a different state $k$ at time $t$ (\ie $k\neq \hat X_{t-}(i)$), then:
\begin{enumerate}
\item if (with probability $1/(n-1)$) we are lucky and $\hat Y_{t-}(i)=k$, choose another unmatched coordinate $j$ uniformly at random, and set $\hat Y_t(j) = \hat X_{t-}(j)$. This decreases $\hat N$ by two; \smallskip
\item if $\hat Y_{t-}(i)\neq k$, set $\hat Y_{t}(i)=k$. This decreases $\hat N$ by one.
\end{enumerate}

\setlength{\fboxsep}{5pt}

\medskip
Now define
    \begin{equation}
        \hat{v}(x,y,t) = \prob{\hat{\tau}>t \,|\, \hat X_0 = x, \, \hat Y_0=y}
    \end{equation}
to be the tail probability of the coupling time $\hat \tau = \hat{\tau}_{n,d}$ under $\hat{c}$. The main result of this paper is the following generalisation of Theorem 3.1 of \cite{Connor.Jacka-2008}.

\begin{thm}\label{thm:main-result}
For any states $x,y\in\hypK$ and time $t\geq 0$,
\begin{equation}
\label{eqn:main-thm}
\hat{v}(x,y,t) = \inf_{c\in \mathcal{C}} \prob{\tau^c>t \,|\, X^c_0=x,\, Y^c_0=y} \,.
\end{equation}
\end{thm}

In other words, $\hat{\tau}$ is the stochastic minimum of all co-adapted coupling times for the pair $(X,Y)$. The proof of this theorem can be found in Appendix~\ref{sec:appendix}.


\section{Limiting behaviour}\label{sec:limits}

Now that we have established the existence of an optimal co-adapted coupling $\hat c$, a natural question to ask is whether or not this coupling is also maximal. (This was answered in the negative when $n=2$ by~\cite{Connor.Jacka-2008}.) Denote by $\pi_n^d$ the uniform distribution on $\hypK$ (recall that $\pi_n^d$ is the equilibrium distribution of $X$), and by $\tau_{n,d}^*$ the maximal coupling time for the pair $(X,Y)$ where $X_0=\mathbf{0}$ and $Y_0\sim\pi_n^d$. The following result is a simple generalisation of Proposition 1 of \cite{Diaconis1990a}:

\begin{lem}\label{lem:tv-dist}
Let 
\[ T_{n,d} = \frac{1}{2}\bra{\frac{n-1}{n}} \log d \,. \]
Then as $d\rinf$, for all $\theta\in\RR$,
\begin{equation}\label{eqn:cutoff}
\norm{\law(X_{T_{n,d}+\theta}) - \pi_n^d}_{\textup{TV}} = 2\Phi\bra{\frac{\sqrt{n-1}}{2}\, e^{-n\theta/(n-1)}}-1+o(1) \,, 
\end{equation}
where $\Phi(\cdot)$ is the standard normal distribution function.
\end{lem}
This shows that the distance between $\law(X)$ (the law of $X$) and $\pi_n^d$  exhibits a \emph{cutoff phenomenon} \cite{Aldous1983,Diaconis-1988,Diaconis-1996} at time $T_{n,d}$, the window size being $O(1)$. Thus $\mathbb{E}[\tau^*_{n,d}]\sim T_{n,d}$.

On the other hand, 
we can bound $\mathbb{E}[\hat{\tau}_{n,d}]$ as follows. As in Appendix~\ref{sec:appendix}, we write $\hat\lambda_t(k,k+s)$ for the rate (according to $\hat Q(t)$) at which $\hat N_\tmi$ jumps from $k$ to $k+s$. Under $\hat{c}$, $\hat{N}$ is a decreasing process, with jumps being of size -1 or -2; the total rate at which $\hat{N}$ jumps is equal (by \C2) to
\begin{equation}\label{eqn:total-rate}
\lc{k,k-2} +\lc{k,k-1} = \frac{k}{n-1}+\frac{k(n-2)}{n-1} = k \,. 
\end{equation}

Now let $M$ be a process that takes only steps of size $k\to k-2$ at rate $k$, and let $\tau_M$ be the time taken for $M$ to be absorbed at zero. If $\hat{N}_0 = M_0 = 2m$ then $\Ex{\tau_M|M_0=2m}\leq \Ex{\hat{\tau}_{n,d}|\hat N=2m}$, thanks to Lemma~\ref{lem:v-increasing} of Appendix~\ref{sec:appendix}. Furthermore,
\[ \Ex{\tau_M|M_0=2m} = \sum_{i=1}^{m} (2i)^{-1} \sim \frac{1}{2}\log m \,. \]

Since $\hat{N}_0=\abs{\hat X_0-\hat Y_0}\sim\text{Bin}(d,(n-1)/n)$, Chebychev's inequality implies that $\hat{N}_0$ (and thus $M_0$) is concentrated around $d(n-1)/n$, and so
\[ \Ex{\tau_M}\sim  \frac{1}{2} \log d  \quad\text{as $d\rinf$} \,. \]

Therefore 
\[ \frac{\Ex{\hat{\tau}_{n,d}}}{\mathbb{E}[\tau^*_{n,d}]} \geq \frac{\Ex{\tau_M}}{\mathbb{E}[\tau^*_{n,d}]} \sim 
\frac{n}{n-1}\,,   \]
and so the optimal co-adapted coupling is not maximal for any fixed $n$.

\medskip
Let us now consider what happens if we let $n\rinf$ while keeping $d$ fixed. Suppose that the $n$ points of $K_n$ are equally spaced on the unit interval $[0,1)$, at locations $\set{0,1/n,2/n,\dots}$. As $n\rinf$ the random walk $X$ on \hypK, with $X_0=\mathbf{0}$, converges in distribution to the random walk $\tilde{X}$ on $[0,1)^d$ for which each coordinate jumps, at incident times of an independent unit-rate Poisson process, to a new location distributed uniformly on $[0,1)$. The equilibrium distribution of $\tilde{X}$ is of course $\pi_\infty^d=\unifdist[0,1)^{\otimes d}$.

\begin{lem}\label{lem:limiting}
For $d$ fixed, as $n\rinf$ the optimal co-adapted coupling $\hat{c}_{n,d}$ of Section~\ref{sec:co-adapted} tends to a maximal coupling.
\end{lem}

\begin{proof}
Let $A_0$ be the set of points in $[0,1)^d$ which have at least one coordinate equal to $0$. Then, by definition of total variation distance,
\begin{align}
\norm{\law(\tilde{X}_t) - \pi_\infty^d}_{\textup{TV}} &= \sup_{A\subset[0,1)^d}\bra{\prob{\tilde{X}_t\in A} - \pi_\infty^d(A)} \nonumber \\
&= \prob{\tilde{X}_t\in A_0} \nonumber \\
&= 1-\prob{\text{all coordinates of $\tilde{X}$ have jumped by time $t$}} \nonumber \\
&= 1-(1-e^{-t})^d \,. \label{eqn:tail_infty}
\end{align}

Now consider the optimal co-adapted coupling strategy $\hat{c}_{n,d}$ as $n\rinf$. Let $\tilde{Y}_0\sim\pi_\infty^d$, and note that $\abs{X_0 - Y_0} = d$ almost surely. From (\ref{eqn:lam2}) we see that $\lambda_t^{\hat{c}_{n,d}}(m,m+s)$, the rate under $\hat{c}_{n,d}$ at which $\hat N_{t-}$ jumps from $m$ to $m+s$, satisfies (when $n\geq 4$)
\begin{align*}
 \lambda_t^{\hat{c}_{n,d}}(m,m+s) &=
\begin{cases}
\frac{m(n-2)}{n-1} & \quad\text{if $s=-1$} \\
\frac{m}{n-1} & \quad\text{if $s=-2$}
\end{cases} \qquad\text{when $m\geq 2$,} \\
 \lambda_t^{\hat{c}_{n,d}}(1,0) &= \frac{n}{n-1} \,,
\end{align*}
with all other rates equal to zero. Let $B$ be a process that takes steps of size $m\to m-1$ at rate $m(n-2)/(n-1)$. Since $\hat{N}_0 = d$, the time taken for $\hat{N}$ to be absorbed at zero (equal to $\hat \tau_{n,d}$) is bounded above by the time taken for $B$ to hit zero from $d$ (again thanks to Lemma~\ref{lem:v-increasing}). But this time, $\tau_B$, is simply equal to the maximum of $d$ independent $\expdist((n-2)/(n-1))$ random variables, and so 
\[ \prob{\tau_B>t \,|\, B_0=d} = 1-\bra{1-e^{-(n-2)t/(n-1)}}^d \,. \]

Finally, using the coupling inequality to lower bound the tail distribution of $\hat\tau_{n,d}$, we see that 
\[ 1-(1-e^{-t})^d \leq \prob{\hat\tau_{n,d}>t} \leq 1-\bra{1-e^{-(n-2)t/(n-1)}}^d \,. \]
As $n\to\infty$ it follows that the coupling time $\hat\tau_{\infty,d}$ achieves equality in the coupling inequality, showing $\hat c_{\infty,d}$ to be a maximal coupling.


\end{proof}

Furthermore, if we now let $d\rinf$, the distance between $\law(\tilde{X}_t)$ and $\pi_\infty^d$ again obeys a cutoff phenomenon, with cutoff time equal to $\log d$. (This may appear surprising, since $T_{n,d} \to \tfrac{1}{2}\log d$ as $n\rinf$ in Lemma~\ref{lem:tv-dist}. However, note that the expression on the right-hand-side of (\ref{eqn:cutoff}) tends to one for all $\theta\in\RR$ as $n\rinf$, showing that $\tfrac{1}{2}\log d$ is not the cutoff time for the limiting process.)

\section{Final comments}\label{sec:comments}

\subsection{A faster coupling for random walks on \hypK}\label{ssec:maximal_on_hypK}
We have demonstrated the existence of (and explicitly described) a maximal co-adapted coupling $\hat c$ for a pair of continuous-time symmetric random walks on \hypK; this has been seen to be submaximal for any fixed $n$, but to converge naturally to a maximal coupling as $n\to\infty$. The reason that $\hat c$ is submaximal can be found in \C2(i), where the strategy for making two new matches is described. Recall that if $\hat X(i)$ jumps to state $k$ at time $t$, where $\hat Y_{t-}=k$, then we choose a different coordinate $j\in U^{\hat c}_t$ and set $\hat Y_t(j) = \hat X_{t-}(j)$. In order for the coupling to be co-adapted, $j$ must be chosen in a way that depends only on the information in $\mathcal{F}_t$. \cite{Matthews-1987} shows how to use information about the future of one process to more efficiently pair coordinates when $n=2$, producing a non-co-adapted (and near-maximal) coupling. This, combined with the result of Lemma~\ref{lem:limiting}, motivates the following suggestion for a faster coupling on \hypK.

\smallskip
Suppose that $X_0 = \mathbf{0}$ and $Y_0 = \mathbf{1}$. (Vectors of zeros and ones respectively -- since $K_n$ is a complete graph, the coupling time is unaffected by coordinate-wise relabelling of the non-zero states of $Y_0$.) Run $X$ until time $\tau = \min\set{\tau_1, \tau_2}$, where
\begin{align*}
\tau_1 &= \inf\set{t: \text{all coordinates have visited the set $\set{2,3,\dots,n-1}$ by time $t$}} \\
\tau_2 &= \inf\set{t: \#\set{i:X_t(i) = 0} = \#\set{i:X_t(i) = 1}}\,,
\end{align*}
and where $\#$ denotes cardinality. Now, if $\tau=\tau_1$, then each coordinate has visited a state not equal to 0 or 1 by time $\tau$: thus $Y$ can be evolved over the period $[0,\tau]$ with every coordinate being independently coupled using the idea depicted in part b) of Figure~\ref{fig:intuition}. If $\tau=\tau_2$ however, then let 
\[ A_\tau^0 = \set{i: X_{\tau}(i)=0}\,, \quad A_\tau^1 = \set{i: X_{\tau}(i)=1}\,, \quad\text{and} \quad A_\tau^+ = \set{i: X_{\tau}(i)\geq 2} \,. \]
On the event $\set{\tau=\tau_2}$, $\#A_\tau^0 = \#A_\tau^1$, and so we may choose some arbitrary bijection $\rho:A_\tau^0\to A_\tau^1$. Now run $Y$ over $[0,\tau]$ as follows: if $i\in A_\tau^0$ then let $Y(i)$ copy the path of $X(\rho(i))$, but with jumps $k\to 0$ replaced by $k\to 1$ and vice versa; if $i\in A_\tau^1$ then let $Y(i)$ similarly copy the path of $X(\rho^{-1}(i))$; coordinates $i\in A_\tau^+$ are coupled independently, as when $\tau=\tau_1$. 

It is simple to check that this produces an adapted coupling of $X$ and $Y$, with coupling time $\tau$. Note that when $n=2$, $\tau = \tau_2$ and we obtain the coupling of \cite{Matthews-1987}. Moreover, $\prob{\tau=\tau_1}\to 1$ as $n\to\infty$, and in the limit we find the maximal coupling $\hat{c}_\infty$ from the proof of Lemma~\ref{lem:limiting}.

\subsection{Random walk on $\mathbb{Z}_n^d$}

The other natural generalisation of the random walk on \hyp is of course to the random walk on $\mathbb{Z}_n^d$, for $n\geq 2$. This is messier to analyse using the stochastic control arguments of this paper, since $\abs{X_t(i)-Y_t(i)}$ can now be greater than 1; although the optimal co-adapted coupling strategy (if one exists) will again be invariant under coordinate permutation, the value function $\hat v$ will now depend not only upon $N_t$, but upon the entire vector $(\zeta_t(0),\zeta_t(1),\dots,\zeta_t(n/2))$, where 
\[ \zeta_t(k) = \#\set{i:\abs{X_t(i)-Y_t(i)} = k}\,. \]

It seems obvious that if $d=1$, the best co-adapted coupling policy is for $Y$ to jump at the same time as $X$, but in the opposite direction around the circle $\mathbb{Z}_d$, with the exception that if $\abs{X_t - Y_t}=1$ then both processes should evolve independently until this equality is broken. When $d\geq 2$, heuristics point towards the best strategy being to couple each coordinate independently (using reflection around the circle), \emph{unless} $\zeta_t(1)$ happens to be even; in this case, we should \emph{pair} those coordinates which are 1 away from coupling, and try to make two new matches at once (in the style of part a) of Figure~\ref{fig:intuition}).

Finally, note that if we appropriately rescale these processes then the coordinates of $X$ converge as $n\rinf$ to a set of $d$ independent Brownian motions on the unit circle. Coupling each coordinate independently by reflection -- the limit of the suggested optimal co-adapted coupling -- should once again be \emph{maximal} for the limiting process.


\bigskip
\appendix
\section{}\label{sec:appendix}

\begin{proof}[Proof of Theorem~\ref{thm:main-result}]

From Definition~\ref{def:Q-hat} it is evident that $\hat{c}$ is invariant under coordinate permutation, and that $\hat{v}(x,y,t)$ only depends on $(x,y)$ through $\abs{x-y}$ (the Hamming distance between $x$ and $y$), and so we shall write
    \[ \hat{v}(m,t) = \prob{\hat{\tau}>t \,|\, \hat N_0 = m} \,, \]
with the convention that $\hat{v}(m,t)=0$ for $m\leq 0$.

 Following the paper~\cite{Connor.Jacka-2008}, we shall write $\lambda^c_t(m,m+s)$ for the rate (according to $Q^c(t)$) at which $N_t^c$ jumps from $m$ to $m+s$, for $s\in\set{-2,\dots,2}$. For example:
\begin{equation}\label{eqn:lambda(m,m-2)}
  (n-1)\lambda_t^c(m,m-2)= \sum_{\substack{i,j\in U^c_{t-} \\i \neq j}}q^c_{(i,Y^c_{t-}(i))(j,X^c_{t-}(j))}(t)
\end{equation}
and
\begin{align}
  (n-1)\lambda_t^c(m,m-1) &= \sum_{i\in U^c_{t-}}\sum_{0\leq k\leq n-1} q^c_{(i,k)(i,k)}(t)\label{eqn:lambda(m,m-1)} \\
  & \qquad + \sum_{\substack{i,j\in U^c_{t-} \\i \neq j}} \, \sum_{0\leq
    l,k\leq n-1} q^c_{(i,k)(j,l)}(t)
  \bra{\indev{k=Y^c_{t-}(i), \,l\neq X^c_{t-}(j)} + \indev{k\neq Y^c_{t-}(i), \,l= X^c_{t-}(j)}} \nonumber  \\
  & \qquad + \sum_{\substack{i\in U^c_{t-} \\j\in M^c_{t-}}} q^c_{(i,Y^c_{t-}(i))(j,Y^c_{t-}(j))}(t)  + \sum_{\substack{i\in M^c_{t-} \\j\in U^c_{t-}}} q^c_{(i,X^c_{t-}(i))(j,X^c_{t-}(j))}(t)  \,.  \nonumber
\end{align}
Similar decompositions may be written down for $\lambda^c_t(m,m+1)$ and $\lambda^c_t(m,m+2)$, but we will have no need of them in the sequel.

The expression for $\lambda_t^c(m,m-2)$ in \eqref{eqn:lambda(m,m-2)} is easy to understand: $N^c$ decreases by two if and only if different unmatched coordinates on $X^c$ and $Y^c$ flip at the same instant, with each flip making one new match. 

$\lambda_t^c(m,m-1)$ comprises four sums however, and so requires a little more explanation. The first sum in \eqref{eqn:lambda(m,m-1)} gives the rate at which the same unmatched coordinate flips on both $X^c$ and $Y^c$ to the same value, making one new match. The second term is the rate at which an unmatched coordinate on one process flips to make a new match, while a different unmatched coordinate on the other process jumps (or possibly stays at its current value) \emph{without} making another new match. Finally, the third and fourth sums in \eqref{eqn:lambda(m,m-1)} give the rate at which an unmatched coordinate on one process flips and makes a new match, while on the other process a matched coordinate is selected and made to stay at its current value.

Using the constraints on the row sums of $Q^c_t$, it is possible to bound the terms in \eqref{eqn:lambda(m,m-2)} and \eqref{eqn:lambda(m,m-1)} as follows:
\begin{align}
(n-1)\lambda_t^c(m,m-2) &= \sum_{\substack{i,j\in U^c_{t-} \\i \neq j}}q^c_{(i,Y^c_{t-}(i))(j,X^c_{t-}(j))}(t)
\leq \abs{U^c_{t-}} = m \,;  \label{eqn:bound-on-2}  \\
\intertext{and similarly,}
 (n-1)\lambda_t^c(m,m-1) &\leq nm. \nonumber
\end{align}

 Moreover, since $Q^c_t$ is doubly stochastic:
\begin{align}
  (n-1)\lambda_t^c(m,m-1) &+ 2(n-1)\lambda_t^c(m,m-2) \nonumber  \\
&\leq \sum_{i\in U^c_{t-}}\bra{\sum_{j=1}^d\sum_{0\leq k,l<n} q^c_{(i,k)(j,l)}(t) \indev{k\neq X^c_{t-}(i),\, k\neq Y^c_{t-}(i)}} \nonumber \\
  & \qquad + \sum_{i\in U^c_{t-}}\bra{\sum_{j=1}^d \sum_{l=0}^{n-1} q^c_{(i,k)(j,l)}(t) \bra{ \indev{j\neq i}+ \indev{j=i, \, l=Y^c_{t-}(i)}}} \nonumber \\
  & \qquad + \sum_{j\in U^c_{t-}}\bra{\sum_{i=1}^d \sum_{k=0}^{n-1} q^c_{(i,k)(j,l)}(t) \bra{ \indev{i\neq j}+ \indev{i=j, \, k=X^c_{t-}(j)}}} \nonumber \\
  & \leq m(n-2) + m + m = nm \,. \nonumber
\end{align}

Denote by $L_n$ the set of nonnegative $\lambda$ satisfying this key linear constraint
\begin{equation} \label{eqn:upper-bound}
(n-1)\lambda(m,m-1) + 2(n-1)\lambda(m,m-2) \leq nm \,.
\end{equation} 
When $n=2$ this reduces to the constraint of \cite{Connor.Jacka-2008}: $\lambda(m,m-1)+2\lambda(m,m-2) \leq 2m$.

\begin{prop}
\label{prop:lambda}
Under $\hat{c}$ the following set of equations hold:
\begin{align}
\hat \lambda_t(m,m+1) &=  \hat \lambda_t(m,m+2) = 0 \,; \label{eqn:lam1} 
\intertext{if $m$ is even, or if $m\geq 2(n-1)/(n-2)$ then}
(n-1)\hat \lambda_t(m,m-2) &= m \quad\text{and} \quad (n-1)\hat \lambda_t(m,m-1) = (n-2)m \,; \label{eqn:lam2} 
\intertext{while if $m$ is odd and $m<2(n-1)/(n-2)$ then}
\hat \lambda_t(m,m-2) &= 0 \quad\text{and} \quad (n-1)\hat \lambda_t(m,m-1) = nm \,. \label{eqn:lam3}
\end{align}
\end{prop}

\begin{proof}
Equation~(\ref{eqn:lam1}) is an immediate consequence of \C1, which implies that no matches are ever broken under $\hat{c}$. When $\hat N_\tmi = m$ is even or satisfies $m\geq 2(n-1)/(n-2)$, it follows from \C2(i) and equation~(\ref{eqn:lambda(m,m-2)}) that 
\[ (n-1)\hat \lambda_t(m,m-2) = \sum_{\substack{i,j\in \hat U_\tmi \\i \neq j}}\hat{q}_{(i,\hat Y_{t-}(i))(j,\hat X_{t-}(j))}(t) = \sum_{\substack{i,j\in \hat U_\tmi \\i \neq j}} \bra{\frac{1}{m-1}} = m \,.  \]
Finally, \C1 and \C2 imply that, under $\hat{c}$, the only non-zero term in equation~(\ref{eqn:lambda(m,m-1)}) is the first sum, and so
\[ (n-1)\hat \lambda_t(m,m-1) =  \sum_{i\in \hat U_\tmi}\sum_{k=0}^{n-1}\hat{q}_{(i,k)(i,k)}(t) \,. \]
Substituting the values of $\hat{q}_{(i,k)(i,k)}(t)$ from \C2 and \C3 completes the proof.
\end{proof}

It follows that the upper bound of (\ref{eqn:upper-bound}) is always attained under $\hat{c}$. Although the framework laid out in Section~\ref{sec:co-adapted} for describing a general coupling $c\in\mathcal{C}$ differs from the setup in \cite{Connor.Jacka-2008}, we can immediately obtain the result of that paper:

\begin{corollary}
\label{cor:d=2}
Theorem~\ref{thm:main-result} holds when $n=2$.
\end{corollary}

\begin{proof}
When $n=2$, Proposition~\ref{prop:lambda} shows that for all $m\in\NN$:
\begin{align*}
\hat \lambda_t(m,m+1) &=  \hat \lambda_t(m,m+2) = 0 \,,\\
\hat \lambda_t(2m,2m-2) &= 2m, \quad\text{and} \quad \hat \lambda_t(2m-1,2m-2) = 2(2m-1) \,.
\end{align*}
The optimality of $\hat{c} = \hat{c}_{2,d}$ now follows from the proof of Theorem 3.1 of \cite{Connor.Jacka-2008}.
\end{proof}

For a strategy $c\in\mathcal{C}$, define the process $S_t^c$ by
    \[ S_t^c = \hat{v}\bra{X_t^c, Y_t^c, T-t}\,, \]
where $T>0$ is some fixed time. This is the conditional probability of $X$ and $Y$ not having coupled by time $T$, when strategy $c$ has been followed over the interval $[0,t]$ and $\hat{c}$ has then been used from time $t$ onwards.  As in \cite{Connor.Jacka-2008}, the optimality of $\hat{c}$ will follow by Bellman's principle \cite{Krylov-1980} if it can be shown that $S^c_{t\wedge\tau^c}$ is a submartingale for all $c\in\mathcal{C}$ (where $s\wedge t=\min\set{s,t}$).


Once again, $S^c$ satisfies
    \begin{equation}\label{eqn:diff-eqn}
        \d S_t^c = \d Z^c_t  + \bra{\mathcal{A}_t^c \hat{v} - \frac{\partial \hat{v}}{\partial t}} \d t \,,
    \end{equation}
where $Z^c_t$ is a martingale, and $\mathcal{A}_t^c$ is the {\lq\lq}generator{\rq\rq} corresponding to the matrix $Q^c(t)$. 
Since $\hat{v}$ is invariant under coordinate permutation, and the Poisson processes $\Lambda_{(i,k)(j,l)}$ are independent, 
\begin{equation*}\label{eqn:generator}
\mathcal{A}_t^c \hat{v}(m,t) =  \sum_{s=-2}^2\lambda_t^c(m,m+s) \squ{\hat{v}(m+s,t) - \hat{v}(m,t)} \,.
\end{equation*}
 To prove that $S^c$ is a submartingale it suffices to show that $\mathcal{A}_t^c\hat{v}$ is minimised by setting $c=\hat{c}$. 
 Thus we seek to maximise over $\lambda\in L_n$, for all $m\geq 0$ and all $t\geq 0$, 
\begin{equation}\label{eqn:maximisation-problem}
\sum_{s=-2}^2 \lambda(m,m+s) \squ{\hat{v}(m,t) - \hat{v}(m+s,t)} \,. 
\end{equation}
Our first step is to simplify this maximisation problem by showing that $\hat{v}(m,t)$ is strictly increasing in $m$. When $n=2$, this result follows trivially from the explicit representation of $\hat{\tau}$ given in \cite{Connor.Jacka-2008}. For $n>2$ however, the result is less obvious and requires a formal proof. 


\begin{lem}
\label{lem:v-increasing}
The tail probability $\hat{v}(m,t)$ is strictly increasing in $m$. 
\end{lem}

\begin{proof}
 We detail here the proof for the case when $n\geq 4$ (for which case \C2 of Definition~\ref{def:Q-hat} applies for all $m>1$): the proof when $n=3$ is similar, using the remark following equation (\ref{eqn:v-recursive}) whenever \C3 applies (\ie when $m\in\set{1,3}$).

We begin by considering $\hat{v}(1,t)$. By (\ref{eqn:lam1}) and (\ref{eqn:lam3}) it follows directly that for all values of $n$,
\begin{equation}
\label{eqn:v(1)}
\hat{v}(1,t) = \exp\bra{-\frac{nt}{n-1}}\,.
\end{equation}

Now consider (for $m>1$) that part of the coupling $\hat{c}$ described in \C2. From (\ref{eqn:lam2}), the total rate at which $\hat N_\tmi$ can change under \C2 is given by 
\[ \lc{m,m-2} +\lc{m,m-1} = \frac{m}{n-1}+\frac{(n-2)m}{n-1} = m \,. \]
Using this, along with (\ref{eqn:lam1}) and (\ref{eqn:lam2}), we obtain for $m>1$:
\begin{align} \hat{v}(m,t) &= e^{-mt} + \int_0^t me^{-mu}\squ{\lc{m,m-2}\hat{v}(m-2,t-u) +\lc{m,m-1}\hat{v}(m-1,t-u)}\d u \nonumber \\
&=  e^{-mt} + \int_0^t \frac{me^{-mu}}{n-1}\squ{m\hat{v}(m-2,t-u) +m(n-2)\hat{v}(m-1,t-u)} \d u \,. \label{eqn:v-recursive}
\end{align}
(A similar expression can be obtained for $\hat{v}(m,t)$ under \C3, noting that the total rate at which $\hat N_\tmi$ can change in this case is $nm/(n-1)$.)

Define $\hat{V}^m(\alpha)$ to be the Laplace transform of $\hat{v}(m,\cdot)$:
\[ \hat{V}^m(\alpha) = \int_0^\infty e^{-\alpha t}\hat{v}(m,t) \d t \,. \]
It then follows from (\ref{eqn:v-recursive}) that, for $m>1$,
\begin{equation}\label{eqn:V-recursive-first}
\hat{V}^m(\alpha) = \frac{1}{m+\alpha} + \frac{1}{(n-1)(m+\alpha)}\bra{m \hat{V}^{m-2}(\alpha) + m(n-2)\hat{V}^{m-1}(\alpha)} \,,
\end{equation}
and so (rearranging)
\begin{equation}
\label{eqn:V-recursive}
(m+\alpha)(n-1)\hat{V}^m(\alpha) = (n-1) + m \hat{V}^{m-2}(\alpha) + m(n-2)\hat{V}^{m-1}(\alpha) \,.
\end{equation}

We need to show that $r(m,t)\geq 0$ for all $m\geq 1$ and $t\geq 0$, where 
\[ r(m,t)= \hat{v}(m,t) - \hat{v}(m-1,t)\,. \]
 By the Bernstein-Widder theorem (see \cite{Feller-1971}, Theorem 1a, Chapter XIII.4), this is equivalent to showing that $R^m(\alpha)$ is totally monotone (TM), where $R^m(\alpha)$ is defined for $\alpha\geq 0$ by 
\[ R^m(\alpha) = \int_0^\infty e^{-\alpha t}r(m,t) \d t \,. \]

 We begin by showing that this is true when $m=1$ and $m=2$, and then use induction. From (\ref{eqn:v(1)}) we see that 
\begin{equation}
\label{eqn:D_alpha(1)}
R^1(\alpha) = \hat{V}^1(\alpha) = \frac{n-1}{n+\alpha(n-1)} \,, 
\end{equation}
and is therefore \tm.

Furthermore, using (\ref{eqn:V-recursive-first}) we obtain
\begin{equation}
\label{eqn:D_alpha(2)}
R^2(\alpha) = \frac{1}{2+\alpha}+ \frac{2(n-2)}{(2+\alpha)(n-1)}\hat{V}^1(\alpha) -  \hat{V}^1(\alpha) =  \frac{n-2}{(2+\alpha)(n+\alpha(n-1))} \,.
\end{equation}
Since the product of two \tm\ functions is itself \tm~\cite{Feller-1971}, it follows that $R^2(\alpha)$ is also \tm, as desired.

Now suppose that we have already shown $R^m(\alpha)$ and $R^{m-1}(\alpha)$ to be \tm, for some $m\geq 2$. Subtracting $(m+\alpha)(n-1)\hat{V}^{m-1}(\alpha)$ from both sides of (\ref{eqn:V-recursive}) yields
\begin{equation}
\label{eqn:v-inc1}
 (m+\alpha)(n-1)R^m(\alpha) = (n-1)\squ{1-\alpha \hat{V}^{m-1}(\alpha)}-mR^{m-1}(\alpha) \,.
\end{equation}
Substituting $m+1$ for $m$ in this expression we obtain
\begin{equation}
\label{eqn:v-inc2}
 (m+1+\alpha)(n-1)R^{m+1}(\alpha) = (n-1)\squ{1-\alpha \hat{V}^m(\alpha)}-(m+1)R^m(\alpha) \,,
\end{equation}
and then subtracting (\ref{eqn:v-inc1}) from (\ref{eqn:v-inc2}) yields
\begin{align}
(m+1+\alpha)(n-1)R^{m+1}(\alpha) &= (n-1)\alpha\squ{\hat{V}^{m-1}(\alpha)-\hat{V}^m(\alpha)}+mR^{m-1}(\alpha) \nonumber \\
&\qquad - \squ{(m+1)-(m+\alpha)(n-1)}R^m(\alpha) \nonumber \\
&= mR^{m-1}(\alpha) + \squ{m(n-2)-1}R^m(\alpha) \,. \label{eqn:v-inc-final}
\end{align}

Since $m\geq 2$ and $n\geq 4$ we see that $m(n-2)-1>0$. Hence, by our induction hypothesis, $R^{m+1}(\alpha)$ can be expressed as the sum of two \tm\ functions, and so is itself \tm. This completes the proof.
\end{proof}

Thus $\hat{v}(m,t) - \hat{v}(m-1,t)\geq 0$ for all $m\geq 0$ and $t\geq 0$. It follows that the terms appearing on the right-hand-side of equation (\ref{eqn:maximisation-problem}) are nonpositive if and only if $s$ is nonnegative. Hence we must set 
\[ \lambda(m,m+1) = \lambda(m,m+2) = 0 \]
in order to achieve the maximum in (\ref{eqn:maximisation-problem}).
It therefore now suffices to maximise 
\begin{align*}
\lambda(m,m-2)&\squ{\hat{v}(m,t)-\hat{v}(m-2,t)} +  \lambda(m,m-1)\squ{\hat{v}(m,t)-\hat{v}(m-1,t)} \\
&= \lambda(m,m-2)\squ{r(m,t)+r(m-1,t)} +  \lambda(m,m-1)r(m,t)
\end{align*}
subject to the constraint from (\ref{eqn:upper-bound}):
\[ (n-1)\lambda(m,m-1) + 2(n-1)\lambda(m,m-2) \leq nm \,. \]
Putting these together we see that we need to maximise (for $m\geq 2$)
\begin{equation}\label{eqn:final-maximisation}
\lambda(m,m-2)\squ{r(m-1,t)-r(m,t)} \,.
\end{equation}

\begin{thm}
\label{thm:main}
For all $n\geq 4$, $m\geq 2$ and $t\geq 0$: $r(m-1,t)\geq r(m,t)$. 
\end{thm}

\begin{proof}
In a similar fashion to the proof of Lemma~\ref{lem:v-increasing}, we show positivity of $r(m-1,t)-r(m,t)$ by showing $R^{m-1}(\alpha) - R^m(\alpha)$ to be \tm, again using induction in $m$. From (\ref{eqn:D_alpha(1)}) and (\ref{eqn:D_alpha(2)}) we see that 
\[ R^1(\alpha) - R^2(\alpha) = \frac{1}{2+\alpha} \]
and so is \tm. Using (\ref{eqn:v-inc-final}) it can be deduced that
\[ R^2(\alpha) - R^3(\alpha) = \frac{n-4}{(2+\alpha)(3+\alpha)(n-1)} \,. \]
Since $n\geq 4$, this difference is also \tm.

Now assume that $R^{m-1}(\alpha) - R^m(\alpha)$ is \tm, for some $m\geq 3$. Substituting $m-1$ for $m$ in (\ref{eqn:v-inc-final}) yields
\begin{equation}
\label{eqn:m->m-1}
(m+\alpha)(n-1)R^m(\alpha) = (m-1)R^{m-2}(\alpha) + \squ{(m-1)(n-2)-1}R^{m-1}(\alpha) \,,
\end{equation}
and subtracting (\ref{eqn:v-inc-final}) from (\ref{eqn:m->m-1}) shows that
\begin{align*}
(m+1+\alpha)(n-1)\squ{R^m(\alpha)-R^{m+1}(\alpha)} &= ((m-1)(n-2)-2)\squ{R^{m-1}(\alpha)-R^m(\alpha)} \\
&\qquad + (m-1)\squ{R^{m-2}(\alpha)-R^{m-1}(\alpha)} \,.
\end{align*}
Finally, since $m\geq 3$ and $n\geq 4$, $(m-1)(n-2)-2\geq 2$ and so it follows from our induction hypothesis that $R^m(\alpha)-R^{m+1}(\alpha)$ is the sum of two \tm\ functions, and hence is itself \tm, as claimed.
\end{proof}

This result allows us to finally complete the proof of Theorem~\ref{thm:main-result} when $n\geq 4$. (The proof for $n=3$ follows a similar line of argument, using an amended version of Theorem~\ref{thm:main}.) When $n\geq 4$, Theorem~\ref{thm:main} and the preceding discussion show that the optimal strategy must maximise $\lambda(m,m-2)$ for $m\geq 2$, and $\lambda(1,0)$. Using the bounds in \eqref{eqn:bound-on-2} and \eqref{eqn:upper-bound}, it follows that this is equivalent to requiring $\lambda(1,0) = n/(n-1)$ and, for $m\geq 2$:
\[ (n-1)\lambda(m,m-2) = m \qquad\text{and}\qquad (n-1)\lambda(m,m-1) = (n-2)m \,. \]
But by Proposition~\ref{prop:lambda}, this is in complete agreement with the rates $\hat \lambda_t$ arising from using our candidate optimal strategy, $\hat{c}$. Thus $\hat{c}$ is truly an optimal co-adapted coupling, as claimed.

\end{proof}


\newpage
\bibliographystyle{apt}
\bibliography{Connor}

\end{document}